\documentclass[12pt, reqno]{amsart}
\usepackage{amsmath, amsthm, amscd, amsfonts, amssymb, graphicx, color}
\usepackage[bookmarksnumbered,plainpages]{hyperref}
\usepackage[all]{xy}

\newtheorem{theorem}{Theorem}[section]
\newtheorem{lemma}[theorem]{Lemma}
\newtheorem{proposition}[theorem]{Proposition}
\newtheorem{corollary}[theorem]{Corollary}
\usepackage{enumerate}
\theoremstyle{definition}

\newtheorem{example}[theorem]{Example}

\newtheorem{Theorem}{\quad Theorem}[section]

\newtheorem{Example}[Theorem]{\quad Example}

\numberwithin{equation}{section}
\input{mathrsfs.sty}
\begin{document}
\title[norm parallelism]{Parallelism in Hilbert $K(\mathcal{H})$-modules}
\author[ M. Mohammadi Gohari and M. Amyari]{M. Mohammadi Gohari and M. Amyari$^*$}
\address{Department of Mathematics, Mashhad Branch,
 Islamic Azad University, Mashhad, Iran.}
  \email{mahdi.mohammadigohari@gmail.com}
\email{maryam\_amyari@yahoo.com  and  amyari@mshdiau.ac.ir}
\subjclass[2010]{Primary 46L08, Secondary 47A30, 47B10, 47B47}
\keywords{Hilbert $C^*$-module, minimal projection, Parallelism, compact operator.\\
*Corresponding author}

\begin{abstract}
Let $(\mathcal{H}, [\cdot, \cdot ])$ be a Hilbert space and $K(\mathcal{H})$ be the $C^*$-algebra of compact operators on $\mathcal{H}$.
In this paper,  we present some characterizations of the norm-parallelism for elements of a  Hilbert $K(\mathcal{H})$-module $\mathcal{E}$ by employing the minimal projections on $\mathcal{H}$. Let $T,S\in \mathcal{L(\mathcal{E})}$. We show that $T \| S$ if and only if there exists a sequence of basic vectors  $\{x_n\}^{\xi_n}$ in $\mathcal{E}$ such that $\lim_n [\langle Tx_n ,  Sx_n \rangle \xi_n , \xi_n ] = \lambda\| T\| \| S\|$
for some $\lambda \in \mathbb{T}$. In addition, we give some equivalence assertions about the norm-parallelism of ``compact'' operators on a Hilbert $C^*$-module.
\end{abstract}
\maketitle
\section{Introduction and preliminaries}
Let $(\mathcal{H}, [\cdot, \cdot ])$ be a Hilbert space, $B(\mathcal{H})$ and $K(\mathcal{H})$
 denote the $C^*$-algebras of all bounded linear operators and compact linear operators on $\mathcal{H}$, respectively.
 A (left) inner product $\mathcal{A}$-module is a left $\mathcal{A}$-module $\mathcal{E}$
  equipped with an $\mathcal{A}$-valued inner product $\langle \cdot , \cdot \rangle$ which is linear
  in the first variable and conjugate linear in the second, and satisfy

  \begin{enumerate}
\item[(i)] $\langle tx,y \rangle=t \langle x,y \rangle,$

\item[(ii)]  $\langle x,y \rangle^*=\langle y,x \rangle,$

\item[(iii)]  $\langle x,x \rangle\geq 0$ and $\langle x,x \rangle=0$ if and only if $x=0,$
\end{enumerate}

 for all $x,y\in \mathcal{E},~~~t \in \mathcal{A}$.
If $\mathcal{E}$ is complete with respect to the norm defined by $\| x\|:=\| \langle x , x \rangle\|^\frac{1}{2} $, then $\mathcal{E}$
is called a Hilbert $\mathcal{A}$-module.
It is called a full Hilbert $\mathcal{A}$-module if the ideal $I={\rm span}\{\langle x,y\rangle:\ x,y\in \mathcal{E}\}$ is dense in $\mathcal{A}$.
 Obviously, every Hilbert space is a Hilbert $\mathbb{C}$-module. Further, every $C^*$-algebra $\mathcal{A}$ can be regarded
 as a Hilbert $C^*$-module over itself via the inner product $\langle a , b \rangle = ab^*$.
For $x \in \mathcal{E}$, we denote by $|x|$ the positive square root of $\langle x, x \rangle$.
  we identify any $C^*$-algebra $\mathcal{A}$ as a $C^*$-subalgebra of $B(\mathcal{H})$ for some Hilbert space $\mathcal{H}$ by using the Gelfand-Naimark theorem. We refer the reader to \cite{La} for
 more details on the theory of Hilbert $C^*$-modules.

 For $\xi, \eta \in \mathcal{H}$, we denote by $\xi \otimes \eta$  the corresponding elementary operator
    $(\xi \otimes \eta)(\nu) = [ \nu, \eta  ] \xi$ for each $\nu \in \mathcal{H}$. Obviously, $\|\xi \otimes \eta\| = \| \xi\|\| \eta\|$.
    A projection $e\in K(\mathcal{H})$ is minimal if $eB(\mathcal{H})e = \mathbb{C}e $. It is easy to show that $e=\xi\otimes \xi$ for some unit vector $\xi \in \mathcal{H}$. The set of all minimal projections is denoted by $P_{m}$. An element $x\in \mathcal{E}$ is said to be a basic vector if $\langle x, x\rangle = e$ for some $e \in P_m$. A system $\{x_{\alpha}\}_{\alpha \in \Lambda}$ of elements  of a Hilbert $K(\mathcal{H})$-module $\mathcal{E}$ is orthonormal if (i)  $\langle x_{\alpha}, x_{\beta} \rangle = 0$ for all $\alpha \neq \beta$,
(ii)  $\langle x_{\alpha} , x_{\alpha} \rangle$ is a minimal projection for each $\alpha$. An orthonormal system is said to be an orthonormal basis for $\mathcal{E}$ if it generates a dense submodule of $\mathcal{E}$. The orthogonal dimension of $\mathcal{E}$ is defined as the cardinal number of any one of its orthonormal bases; see \cite{Bak}. In the rest of the paper, $\{x_n\}^{\xi_n}$  denotes a sequence of basic vectors such that $\langle x_n, x_n \rangle = e_n = \xi_n\otimes \xi_n$, where $\xi_n \in \mathcal{H}$, $\|\xi_n\| = 1$ and $e_n \in P_m$, for all $n\in \mathbb{N}$.

Let $\mathcal{E}$ be a Hilbert $K(\mathcal{H})$-module. A mapping  $T: \mathcal{E} \to \mathcal{E} $ is called adjointable if $\langle Tx, y \rangle = \langle x, T^*y \rangle$ for all $x, y\in \mathcal{E}$. Such a map $T$ is a $K(\mathcal{H})$-linear ($T(ax) = aTx$ for all $a \in K(\mathcal{H})$) and bounded. The set of all adjointable mappings on $\mathcal{E}$ is denoted by $\mathcal{L(\mathcal{E})}$. We denote by $\mathcal{K(\mathcal{E})}$ the closed linear subspace of $\mathcal{E}$ spanned by $\{\theta_{x,y}: x, y \in \mathcal{E}\}$,  where $\theta_{x,y}(z):= \langle z, y \rangle x$. The elements of $\mathcal{K(\mathcal{E})}$ are called ``compact'' operators. It is easy to see that $\theta_{x,y} = \theta_{y,x}^*$ and $\|\theta_{x,x}\| = \|x\|^2$
for all $x, y\in \mathcal{E}$.

Let $(X , \| \cdot\|)$ be a complex normed space and $x , y \in X$. We say that $x$ is orthogonal to $y$ in the Birkhoff--James sense, denoted by $x \perp_{B} y$, if
$$\| x\| \leq \| x + \lambda y\|~~~ \text {for~~~ all }~~~ \lambda \in\mathbb{C}.$$
For more  properties of the Birkhoff--James orthogonality in Hilbert $C^*$-modules; see \cite{AR1, C, C1, G, Gr}.

Let $x , y \in X$. We say that $x$ is norm-parallel to $y$, denoted by  $x \| y$ if
$$\| x + \lambda y\| = \| x\| + \| y\| ~~~ \text{for ~~~some}~~~ \lambda \in \mathbb{T} = \{\alpha \in \mathbb{C}: | \alpha| = 1 \}$$
For more  details one can see \cite{Am,Z1}.
 Recall that the norm-parallelism is a symmetric relation
 and $\mathbb{R}$-homogenous, that is $x \|y \Leftrightarrow \alpha x \| \beta y$ for all $\alpha, \beta \in \mathbb{R}$,
 but not transitive, that is $x \|y$ and $y \| z$ does not imply $x \| z$ in general; see \cite[Example 2.7]{Z2}. In the setting of inner product spaces, norm parallelism is equivalent to linear dependence. In the framework of normed linear spaces  linear dependence of two elements implies the norm parallelism, but the converse is not true in general.
In the rest of paper, $\mathcal{E}$ denotes a Hilbert $K(\mathcal{H})$-module.
 Let $\mathcal{F}$ be a Hilbert $C^*$-module over an arbitrary $C^*$-algebra $\mathcal{A}$  and $x, y \in \mathcal{F}$ . Zamani and Moslehian \cite[Theorem 2.3]{Z2} proved that $x \| y$ if and only if there exists a sequence of unit vectors ${\xi_{n}}$
 in a Hilbert space such that
\begin{align*}
\lim_{n \to \infty}|[ \langle x , y\rangle \xi_{n} , \xi_{n} ]|  = \left\| x \right\|\left\| y \right\|.
\end{align*}
We show that for two elements $x, y \in \mathcal{E}$, $x \| y$  if and only if there exists a unit vector $\xi \in \mathcal{H}$ such that
\begin{align*}
|[ \langle x , y\rangle \xi , \xi ]|  =  \| x \|\| y \|.
\end{align*}

 The paper is organized as follows:\\
 In Section 2, we obtain some characterizations of the norm-parallelism for elements of a Hilbert $K(\mathcal{H})$-module.
 we also prove that linear dependence is equivalent to the norm parallelism for certain Hilbert $K(\mathcal{H})$-module.

In Section 3, we show that for elements $T, S \in \mathcal{L}(\mathcal{E})$, the relation $T \| S$ holds if and only if there exists a sequence of basic vectors  $\{x_n\}^{\xi_n}$ in $\mathcal{E}$ such that
\begin{align*}
\lim_n |[\langle Tx_n ,  Sx_n \rangle \xi_n , \xi_n ]| = \| T\| \| S\|.
\end{align*}
Let $K(\mathcal{H})$ is a finite dimensional $C^*$-algebra and $\mathcal{E}$ be a full Hilbert $K(\mathcal{H})$-module and $T, S \in \mathcal{L}(\mathcal{E})$. Then we show  that $T \| S$ if and only if there exists a basic vector $\{x\}^{\xi}$ in $\mathcal{E}$ such that $|[\langle Tx,  Sx\rangle \xi , \xi ]| = \| T\| \| S\|.$
In addition, we give some equivalence assertions about the norm-parallelism of ``compact'' operators on a Hilbert $C^*$-module.

\section{Norm-parallelism in Hilbert $K(\mathcal{H})$-modules.}

In this section, we obtain some characterizations of the norm-parallelism for elements of a Hilbert $K(\mathcal{H})$-module $\mathcal{E}$.
 Let $ x, y\in \mathcal{E}$.  If there exists a unit vector $\xi \in \mathcal{H}$ such that $|[ \langle x , y \rangle \xi , \xi ]| =\| x\| \| y\|$,  then \cite[Theorem 2.3]{Z2} implies that $x \|y$. The question is under which conditions the converse is true. The following theorem shows that if  $ x, y\in \mathcal{E}$ and $x \|y$, then  there exists a unit vector $\xi \in \mathcal{H}$ such that $|[ \langle x , y \rangle \xi , \xi ]| =\| x\| \| y\|$.

\begin{theorem}\label{T1}
Let $x , y \in \mathcal{E}$. Then the following statements are equivalent:
\begin{enumerate}
\item[(i)] $x \| y$,
\item[(ii)] There exists a unit vector $\xi \in \mathcal{H}$ such that
\begin{align*}
\Big| [\langle x , y \rangle \xi , \xi ]\Big| =\| x\| \| y\|,
\end{align*}
\end{enumerate}
or equivalently
$$\Big| [\langle x , y \rangle \xi , \xi ]\Big| = [\langle x , x\rangle \xi , \xi ]^{\frac{1}{2}}[\langle y , y \rangle \xi , \xi ]^{\frac{1}{2}}.$$
\begin{proof}
$(i) \Rightarrow(ii)$: Let $x \| y$. Thus $\| x + \lambda y\| = \| x\| + \| y\| $ for some $\lambda \in \mathbb{T}$. By \cite[Proposition 2 part (c)]{Bak}, there exists a minimal projection $e \in P_m$ such that $\| x + \lambda y\| = \| e (x + \lambda y)\|$. Hence
\begin{align*}
\| e \langle x + \lambda y, x + \lambda y \rangle e\| = \| e (x + \lambda y)\|^2 = \| x + \lambda y\|^2.
\end{align*}
Since $e = \xi \otimes \xi$ for some unit vector $\xi \in \mathcal{H}$, we obtain
\begin{align*}
[\langle x + \lambda y , x + \lambda y \rangle \xi , \xi ]= (\| x\| + \| y\|)^2,
\end{align*}
whence
\begin{align*}
(\| x\| + \| y\|)^2 = [\langle x + \lambda y , x + \lambda y \rangle \xi , \xi ]
 &=  [\langle x , x \rangle \xi , \xi ] +  2{\rm Re} [\langle x , \lambda y \rangle \xi , \xi ] + [\langle y , y \rangle \xi , \xi ] \\
&\leq  [\langle x , x \rangle \xi , \xi ]+  2\Big| [\langle x , \lambda y \rangle \xi , \xi ]\Big| + [\langle y , y \rangle \xi , \xi ] \\
 &\leq  \| x\|^2 + 2\| x\|\| y\| + \| y\|^2.
\end{align*}
Therefore
$[\langle x , x \rangle \xi , \xi ] = \| x\|^2,~~~ \Big| [\langle x ,  y \rangle \xi , \xi ]\Big| =\| x\| \| y\|$ and $[\langle y , y \rangle \xi , \xi ] = \| y\|^2$.

Hence $\|ex\|=\|x\|,~~~\|ey\|=\|y\|,~~~\|e\langle x , y \rangle e\|=\|\langle x , y \rangle\|$ for some $e \in P_m$.

$(ii) \Rightarrow(i)$: It is evident, see \cite[Theorem 2.3]{Z2}.
\end{proof}
\end{theorem}

Note that if $\mathcal{E}=K(\mathcal{H})$, then \cite[Theorem 2.10]{Z3} is a corollary of the above theorem.

Now we  peresent another characterization of norm-parallelism for elements of $\mathcal{E}$.

\begin{proposition}\label{p1}
Let $x , y \in \mathcal{E}$. Then the following statements are equivalent:
\begin{enumerate}
\item[(i)] $x \| y$,
\item[(ii)] $\lambda\| x\|\| y\|$ is an eigenvalue of $\langle x , y \rangle$ for some $\lambda \in \mathbb{T}$.
\end{enumerate}
\begin{proof}
$(i) \Rightarrow(ii)$: Let $x \| y$. By Theorem \ref{T1} , there exists a unit vectors $\xi \in \mathcal{H}$  such that $|[\langle x , y \rangle \xi , \xi ]| =\| x\| \| y\|$.
 Hence  $\| x\| \| y\|=|[\langle x , y \rangle \xi , \xi ]| \leq \| \langle x , y \rangle \xi\| \leq \| \langle x , y \rangle \| \leq \| x\| \| y\|$, which ensures
$|[\langle x , y \rangle \xi , \xi ]|= \| \langle x , y \rangle \xi\|$. So $ \langle x , y \rangle \xi = \alpha \xi$ for some $\alpha \in \mathbb{C}$,
whence $|\alpha| = \| \langle x , y \rangle \xi\| = \| x\| \| y\| $. Thus $\alpha= \lambda \| x\| \| y\|$ for some $\lambda \in \mathbb{T}$.
Therefore $\langle x , y \rangle \xi = \lambda \| x\| \| y\| \xi$, or equivalently, $\lambda\| x\|\| y\|$ is an eigenvalue of $\langle x , y \rangle$.

$(ii) \Rightarrow(i)$: Let $\lambda\| x\|\| y\|$ be an eigenvalue of $\langle x , y \rangle$ for some $\lambda \in \mathbb{T}$. Hence $\langle x , y \rangle \xi = \lambda \| x\| \| y\| \xi$ for some unit eigenvector $\xi \in \mathcal{H}$. Thus $[\langle x , y \rangle \xi - \lambda \| x\| \| y\| \xi , \xi ]  = 0$ and so $|[\langle x , y \rangle \xi , \xi ]| =\| x\| \| y\|$. It follows from that $(ii) \Rightarrow(i)$ of Theorem \ref{T1} $x \| y$.
\end{proof}
\end{proposition}
\begin{corollary}\label{c2}
Let $x , y \in \mathcal{E}$. Then the following statements are equivalent:
\begin{enumerate}
\item[(i)] $x \| y$,
\item[(ii)] $\langle x , y \rangle \| I$,
\item[(iii)] $|\langle x , y \rangle| \| I$.
\end{enumerate}
\begin{proof}
Take $t=\langle x , y \rangle$. From Proposition \ref{p1} we conclude that $x\|y$ if and only if $\lambda\|t\|$ is an eigenvalue of $t$ for some $\lambda \in \mathbb{T}$.
 The result now follows immediately from  \cite[Theorem 2.21]{Z2}.
\end{proof}

\end{corollary}
The following result is a combination of Proposition \ref{p1} and Theorem 3.10 of \cite{Z1}.
\begin{corollary}\label{c3}
Let $x , y \in \mathcal{E}$. Then the following statements are equivalent:
\begin{enumerate}
\item[(i)] $x \| y$,
\item[(ii)] $\langle x , y \rangle \| I$,
\item[(iii)] $\langle x , y \rangle \| \langle y , x \rangle$,
\item[(iv)] $\langle x , y \rangle^m \| I \,\quad (m \in \mathbb{N}),$
\item[(v)] $\langle x , y \rangle^m \| \langle y , x \rangle^m \,\quad (m \in \mathbb{N}).$
\end{enumerate}
\end{corollary}
As a consequence of Corollary \ref{c3}, we obtain the following result.

\begin{corollary}
Let $x , y \in \mathcal{E}$. Then the following statements are equivalent:
\begin{enumerate}
\item[(i)] $x \| y$,
\item[(ii)] there exists a ``compact'' operator $A \in K(\mathcal{E})$ such that $Ax \| A^*y$,
\item[(iii)] there exists a ``compact" operator $A \in K(\mathcal{E})$ such that $A^mx \| A{^*}^my$ for each ($m \in \mathbb{N})$.

\end{enumerate}
\begin{proof}
$(i)\Leftrightarrow(ii)$: By $(i)\Leftrightarrow(ii)$ of Corollary \ref{c3},
$$ x\| y  \Longleftrightarrow \langle x , y \rangle \| I .$$

Further, by $(ii)\Leftrightarrow(iii)$ and $(iii)\Leftrightarrow(v)$ of Corollary \ref{c3}, we get
\begin{align*}
\langle x , y \rangle \| I &\Longleftrightarrow \langle x , y \rangle \| \langle x , y \rangle^* \Longleftrightarrow \langle x , y \rangle \| \langle y , x \rangle  \\ &\Longleftrightarrow \langle x , y \rangle \langle x , y \rangle \langle x , y \rangle \| \langle y , x \rangle\langle y , x \rangle\langle y , x \rangle \\ &\Longleftrightarrow  \langle \langle x , y \rangle x , \langle y , x \rangle y \rangle \| \langle \langle y , x \rangle y  ,  \langle x , y \rangle x \rangle \\& \Longleftrightarrow \langle \langle x , y \rangle x , \langle y , x \rangle y \rangle \| I.
\end{align*}
Thus, we reach $x \| y \Longleftrightarrow \langle x , y \rangle x \| \langle y , x \rangle y$.
Putting $A =\theta_{x,y}$. Hence  $x \| y \Longleftrightarrow Ax \| A^*y$.

$(i)\Leftrightarrow(iii)$: Utilizing Corollary \ref{c3}, we reach
\begin{align*}
x\| y \Longleftrightarrow \langle x , y \rangle \| I \Longleftrightarrow \langle x , y \rangle \| \langle y , x \rangle &\Longleftrightarrow \langle x , y \rangle^{2m +1} \| \langle y , x \rangle^{2m +1}\\ &\Longleftrightarrow \langle x , y \rangle^m \langle x , y \rangle \langle x , y \rangle^m \| \langle y , x \rangle^m\langle y , x \rangle\langle y , x \rangle^m\\
&\Longleftrightarrow \langle \langle x , y \rangle^m x , \langle y , x \rangle^m y \rangle \| \langle \langle y , x \rangle^m y  ,  \langle x , y \rangle^m x \rangle\\
&\Longleftrightarrow  \langle x , y \rangle^m x \| \langle y , x \rangle^m y \Longleftrightarrow A^mx \| A{^*}^my
\end{align*}
for all $m\in \mathbb{N}$, where $A =\theta_{x,y}$.
\end{proof}
\end{corollary}
\begin{lemma}\cite[Corollary 2.5]{Z2}\label{L1}
Let $x, y \in \mathcal{E}$. If $x \| y$, then there exists $\lambda \in \mathbb{T}$ such that

 $x \perp_{B} (\| y\| \langle x , x \rangle x + \lambda\| x\| \langle y , x \rangle x)$   and   $y \perp_{B} (\|x\| \langle y , y \rangle y + \lambda\| y\|    \langle y , x \rangle y)$.
\end{lemma}
Recall that, in a normed linear space if two elements are linear dependence, then they are norm parallelism, but  in general the converse is not true. In the following example we see that the converse is true for certain Hilbert $K(\mathcal{H})$-modules.
\begin{example}
Suppose that $\mathcal{H}$ is an infinite dimensional Hilbert space. Then $\mathcal{H}$ equipped with
 $\langle x , y \rangle = x \otimes y $ is a left Hilbert $K(\mathcal{H})$-module and its norm
 $\|x\|^2=\|\langle x , x \rangle\|=\|x \otimes x\|$ coincides with the original norm on $\mathcal{H}$.
We can assume that each unit vector $y \in \mathcal{H}$ gives an orthonormal basis for $\mathcal{H}$, since $\langle x , y \rangle y =( x \otimes y)(y) =x$ for each $x \in \mathcal{H}$. Therefore the dimension of $\mathcal{H}$ as a Hilbert module is 1. Let $x , y \in \mathcal{H}$ and $x \neq 0$. If $x \| y$, then  by Lemma \ref{L1}, we have $x \perp_{B} (\| y\| \langle x , x \rangle x + \mu\| x\| \langle y , x \rangle x)$. On the other hand, $\langle x , x \rangle K(\mathcal{H})\langle x , x \rangle = \mathbb{C} \langle x , x \rangle$. It follows from \cite[Proposition 2.4]{AR1} that
\begin{align*}
\langle x , x \rangle\langle  (\| y\| \langle x , x \rangle x + \lambda   \| x\|\langle y , x \rangle x) , x \rangle = 0,
\end{align*}
or equivalently,
\begin{align*}
0 &= (x \otimes x)\bigg((\| y\| (x \otimes x) (x) + \lambda  \| x\| (y \otimes x) (x))\otimes x \bigg) =  (x \otimes x)\bigg((\| y\| \| x\|^2 x + \lambda  \| x\|^3 y)\otimes x \bigg) \\ &= (x \otimes x)(\| y\| \| x\|^2 x \otimes x + \lambda  \| x\|^3 y \otimes x).
\end{align*}
Thus
\begin{align}\label{0}
\| x\|^4\| y\| x \otimes x +  \lambda [ y ,  x ] \|x\|^3 x \otimes x =0.
\end{align}
Taking norm in (\ref{0}), we reach $|[ x ,  y ]| = \|x\|\|y\|$. Thus $x$ and $y$ are linear dependence.
\end{example}

\section{Norm-parallelism of adjointable maps.}
In this section, we investigate some characterizations of the norm-parallelism for pairs of operators on a Hilbert $K(\mathcal{H})$-module $\mathcal{E}$,
we also obtain the characterizations of the norm-parallelism for pairs of operators in $\mathcal{K(\mathcal{E})}$.
The following result was proved  by Baki\'c and Gulja\u s \cite[Lemma 5]{Bak}.
\begin{lemma}\label{L2}
Let $T \in \mathcal{L}(\mathcal{E})$. Then
$\|T\| = \sup\{\|Tx\|: x \in \mathcal{E} ~~~ \text{and}~~~\langle x ,  x \rangle = e \in P_m \}.$
\end{lemma}
Now, by using the above lemma, we give a characterization of the norm-parallelism in $\mathcal{L}(\mathcal{E})$,
\begin{theorem}\label{T2}
Let  $T, S \in \mathcal{L}(\mathcal{E})$. Then the following statements are equivalent:
\begin{enumerate}
\item[(i)] $T \| S$,
\item[(ii)] there exists a sequence of basic vectors  $\{x_n\}^{\xi_n}$ in $\mathcal{E}$ such that
\begin{align*}
\lim_n [\langle Tx_n ,  Sx_n \rangle \xi_n , \xi_n ] = \lambda\| T\| \| S\|
\end{align*}
for some $\lambda \in \mathbb{T}$.
\end{enumerate}
\begin{proof}
$(i) \Rightarrow(ii)$: Let $T \| S$. Then $\| T + \lambda S\| = \|T\| + \| S\| $ for some $\lambda \in \mathbb{T}$. From Lemma \ref{L2} we have

$$\| T + \lambda S\| = \sup \{ \| (T +\lambda S)(x)\|: x \in \mathcal{E}~~~and~~~\langle x ,  x  \rangle = e  \in P_m\}.$$

Hence there exists a sequence of basic vectors  $\{x_n\}^{\xi_n}$ in $\mathcal{E}$ such that
\begin{align*}
 (\|T\| + \| S\|)^2&=\| T + \lambda S\|^2\\ &= \lim_n \| (T + \lambda S)(x_n)\|^2 = \lim_n \| \langle (T + \lambda S)(x_n) , (T + \lambda S)(x_n)\rangle\| \\
 &= \lim_n \| \langle (T + \lambda S)(e_n x_n) , (T + \lambda S)(e_n x_n)\rangle\| \\
 &= \lim_n \| \langle e_n(T + \lambda S)( x_n) , e_n(T + \lambda S)( x_n)\rangle\| \\
&= \lim_n \| e_n\langle (T + \lambda S)(x_n) , (T + \lambda S)(x_n)\rangle e_n\| \\
&= \lim_n [\langle Tx_n + \lambda Sx_n , Tx_n + \lambda Sx_n\rangle \xi_n , \xi_n ]\\
&= \lim_n [\langle Tx_n ,  Tx_n\rangle \xi_n , \xi_n ] + 2  \lim_n {\rm Re}[\langle Tx_n , \lambda Sx_n \rangle \xi_n , \xi_n ] + \lim_n [\langle Sx_n ,  Sx_n \rangle \xi_n , \xi_n ]\\ &\leq  \| T\|^2 + 2 \lim_n |[\langle Tx_n ,  Sx_n \rangle \xi_n , \xi_n ]| + \| S\|^2 \\ &\leq \| T\|^2 + 2 \lim_n \| Tx_n\| \| Sx_n\| + \| S\|^2 \\ &\leq \| T\|^2 + 2 \| T\| \| S\| + \| S\|^2= (\|T\| + \| S\|)^2.
\end{align*}
Thus we get $\displaystyle{\lim_n} |[\langle Tx_n ,  Sx_n \rangle \xi_n , \xi_n ]| = \| T\| \| S\|$, $\displaystyle{\lim_n} [\langle Tx_n ,  Tx_n \rangle \xi_n , \xi_n ] = \| T\|^2$ and \\
 $\displaystyle{\lim_n} [\langle Sx_n ,  Sx_n \rangle \xi_n , \xi_n ] = \| S\|^2$.

$(ii) \Rightarrow(i)$: Suppose that $(ii)$ holds. By the Cauchy-Schwarz inequality we have
\begin{align*}
\lim_n |[\langle Tx_n ,  Sx_n\rangle \xi_n , \xi_n ]|^2 &\leq \lim_n\| \langle Tx_n ,  Sx_n \rangle \xi_n\|^2 \\
 &= \lim_n[\langle Tx_n ,  Sx_n\rangle^*\langle Tx_n ,  Sx_n \rangle  \xi_n , \xi_n ]
  \\ &\leq \| T\|^2\lim_n |[\langle Sx_n ,  Sx_n \rangle \xi_n , \xi_n]|
 \\ &\leq  \| T\|^2 \| S\|^2,
\end{align*}
whence $\displaystyle{\lim_n} [\langle Tx_n ,  Tx_n \rangle \xi_n , \xi_n ] = \| T\|^2$. A similar computation shows that\\
$\displaystyle{\lim_n} [\langle Sx_n ,  Sx_n \rangle \xi_n , \xi_n ] = \| S\|^2$. Hence
\begin{align*}
\| T\|^2 + 2 \| T\| \| S\| + \| S\|^2 &=  \lim_n [\langle Tx_n ,  Tx_n \rangle \xi_n , \xi_n ] + 2 \lim_n {\rm Re}[\langle Tx_n , \lambda Sx_n \rangle \xi_n , \xi_n ] + \lim_n [\langle Sx_n ,  Sx_n \rangle \xi_n , \xi_n ]\\
&= \lim_n [\langle Tx_n +\lambda Sx_n , Tx_n + \lambda Sx_n \rangle \xi_n , \xi_n ]\\
 &\leq \| T + \lambda S\|^2 \leq \| T\|^2 + 2 \| T\| \| S\| + \| S\|^2.
\end{align*}
Thus $\| T + \lambda S\| = \|T\| + \| S\| $ that is $T \| S$.
\end{proof}
\end{theorem}
As an immediate consequence of Theorem \ref{T2} we obtain the following result.
\begin{corollary}\label{Tete}
Let $x, y, z, u \in \mathcal{E}$. Then the following statements are equivalent:
\begin{enumerate}
\item[(i)] $\theta_{x, y} \| \theta_{z, u}$,
\item[(ii)]  there exists a sequence of basic vectors  $\{x_n\}^{\xi_n}$ in $\mathcal{E}$
 such that
\begin{align*}
\lim_n [\langle x_n, y \rangle \langle x, z \rangle \langle u, x_n \rangle \xi_n , \xi_n ] = \lambda \| \theta_{x, y}\| \| \theta_{z, u}\|,
\end{align*}
for some $\lambda \in \mathbb{T}$.
\end{enumerate}
\begin{proof}
 $(i)\Leftrightarrow(ii)$: Applying Theorem \ref{T2} to maps $\theta_{x, y}$ and $\theta_{z, t},$  we obtain the result.
\end{proof}
\end{corollary}
\begin{proposition}\label{p2}
Let $T, S \in \mathcal{L}(\mathcal{E})$. Then the following statements are equivalent:
\begin{enumerate}
\item[(i)] $T \| S$,
\item[(ii)] $T^*T \| T^*S$ and $\| T^*S\| = \| T\|\| S\|$.
\end{enumerate}
\begin{proof}
Suppose that $T \| S$. Then  $\| T + \lambda S\| = \|T\| + \| S\|$ for some $\lambda \in \mathbb{T}$. Let $\pi: \mathcal{L}(\mathcal{E}) \to B(\mathcal{H})$ be a faithful representation of the $C^*$ algebra $\mathcal{L}(\mathcal{E})$. Using \cite[Theorem 3.3]{Z1} we get
\begin{align*}
T \| S &\Longleftrightarrow \pi(T) \| \pi(S) \\
&\Longleftrightarrow (\pi(T)^* \pi(T) \|  \pi(T)^*\pi(S) ~~~ \text{and}~~~ \| \pi(T)^* \pi(S)\| = \| \pi(T)\|\| \pi(S)\|) \\
&\Longleftrightarrow (T^*T \| T^*S \quad  \text{and} \quad  \| T^*S\| = \| T\|\| S\|).
\end{align*}
\end{proof}
\end{proposition}
\begin{corollary}\label{c4}
Let $x , y \in \mathcal{E}$. Then the following statements are equivalent:
\begin{enumerate}
\item[(i)] $T \| I$,
\item[(ii)]  $T \| T^*$
\item[(iii)] There exists a sequence of basic vectors  $\{x_n\}^{\xi_n}$ in $\mathcal{E}$ such that
\begin{align*}
\lim_{n} \| Tx_n - \lambda  \| T\| x_n \| =0.
\end{align*}
\end{enumerate}
\begin{proof}
$(i)\Leftrightarrow(ii)$: It follows from $(i) \Rightarrow(ii)$ of \cite[Theorem 3.10]{Z1} and by utilizing the same reasoning in the proof of  Proposition \ref{p2}.

$(i) \Rightarrow(iii)$: Let $T \| I$. Then there exists a sequence of basic vectors  $\{x_n\}^{\xi_n},$ in $\mathcal{E}$ such that
\begin{align*}
\lim_n [\langle Tx_n ,  x_n \rangle \xi_n , \xi_n ] = \lambda \| T\| \quad \text{and} \quad  \lim_n [\langle Tx_n ,  Tx_n \rangle \xi_n , \xi_n ] = \| T\|^2,
\end{align*}
for some $\lambda \in \mathbb{T}$. Therefore
\begin{align*}
\lim_n \| Tx_n - \lambda  \| T\| x_n \|^2 &= \lim_n \| \langle Tx_n - \lambda  \| T\| x_n , Tx_n - \lambda  \| T\| x_n  \rangle\|\\
 &= \lim_n \| e_n \langle Tx_n - \lambda  \| T\| x_n , Tx_n - \lambda  \| T\| x_n  \rangle e_n \| \\
  &= \lim_n [\langle Tx_n - \lambda  \| T\| x_n , Tx_n - \lambda  \| T\| x_n  \rangle \xi_n , \xi_n ] \\
  &\leq \| T\|^2 - 2|\lambda|^2\| T\|^2+ \| T\|^2= 0.
\end{align*}
 Conversely, suppose that (iii) holds. We have
\begin{align*}
1 + \| T\| \geqslant \| T + \lambda I\| &\geqslant \| Tx_n + \lambda x_n\| \geqslant \| \lambda x_n + \lambda \| T\| x_n - (-Tx_n +  \lambda \| T\| x_n)\| \\&\geqslant  \| \lambda x_n + \lambda \| T\| x_n\| - \| Tx_n -  \lambda \| T\| x_n\| = 1+ \| T\| - \| Tx_n -  \lambda \| T\| x_n\|.
\end{align*}
Taking limits, we obtain the desired result.
\end{proof}
\end{corollary}

Note that  if $ T, S\in \mathcal{L}(\mathcal{E})$ and there exists a basic vector  $x$ in $\mathcal{E}$ such that $[\langle Tx ,  Sx \rangle \xi, \xi] = \lambda\| T\| \| S\|$
for some $\lambda \in \mathbb{T}$. Then $(ii) \Rightarrow(i)$ of Theorem \ref{T2} ensures that $T \| S$. The question is under which conditions the converse is true. To answer this question, we need the following lemma.
\begin{lemma}\cite[Theorem 2.5]{AR2}\label{L3}
Let $\mathcal{F}$ be a full Hilbert $C^*$-module over an arbitrary $C^*$-algebra $\mathcal{A}$. For each norm-bounded sequence $\{x_n\}$ in $\mathcal{F}$, there exists a subsequence $\{x_{n_k}\}$ of $\{x_n\}$ and a element $x\in \mathcal{F}$ such that
 $$\lim_n \| \langle x_{n_k}, y\rangle - \langle x, y\rangle\| = 0 \quad \text{for all}\quad y \in \mathcal{F},$$
if and only if $\mathcal{A}$ is a finite dimensional $C^*$ -algebra.
\end{lemma}
The next result reads as follows.

\begin{proposition}\label{P1}
Let  $K(\mathcal{H})$ be a finite dimensional $C^*$-algebra and $\mathcal{E}$ be a full Hilbert $K(\mathcal{H})$-module
and $T, S \in \mathcal{L}(\mathcal{E})$. Then the following statements are equivalent:
\begin{enumerate}
\item[(i)] $T \| S$,
\item[(ii)] there exist a basic vector $x$ in $\mathcal{E}$  and a unit vector $\xi \in  \mathcal{H}$ such that
\begin{align*}
[\langle Tx ,  Sx \rangle \xi, \xi] = \lambda\| T\| \| S\|
\end{align*}
for some $\lambda \in \mathbb{T}$,
\item[(iii)] there exists a basic vector  $x$ in $\mathcal{E}$ such that $\| Tx\| = \| T\|$, $\| Sx\| = \| S\|$ and $Tx \| Sx$.
\end{enumerate}
\begin{proof}
$(i) \Rightarrow(ii)$: Let $T \| S$. Then  $\| T + \lambda S\| = \|T\| + \| S\| $ for some $\lambda \in \mathbb{T}$.
 By the proof of implication $(i) \Rightarrow(ii)$ of Theorem \ref{T2}, there exists a sequence of basic vectors  $\{x_n\}^{\xi_n}$ in $\mathcal{E}$ such that
\begin{align*}
\| T + \lambda S\| = \lim_n \| (T + \lambda S)(x_n)\|.
\end{align*}

Since $\{x_n\}$ is a bounded sequence and  the $C^*$-algebra $K(\mathcal{H})$ is finite dimensional, by Lemma \ref{L3}, there is a subsequence $\{x_{n_k}\}$ of $\{x_n\}$ which converges to a unit vector $x$  in $\mathcal{E}$.
\begin{align*}
\lim_n \| \langle x_{n_k}, y\rangle - \langle x , y\rangle\| = 0,
\end{align*}
for all $y \in \mathcal{E}$. Hence
\begin{align*}
\lim_n \| \langle (T + \lambda S)(x_{n_k}), y\rangle - \langle (T + \lambda S)(x), y\rangle\| = 0.
\end{align*}
Thus
\begin{align*}
\| T + \lambda S\| = \lim_n \| (T + \lambda S)(x_n)\| = \| (T + \lambda S)(x)\|.
\end{align*}
Using a similar argument as in the proof of implication $(i) \Rightarrow(ii)$ of Theorem \ref{T2} we obtain the desired result.

$(ii) \Rightarrow(iii)$: Suppose that (ii) holds. Then
\begin{align*}
\| T\| \| S\| = |[\langle Tx ,  Sx \rangle \xi, \xi]| \leq \| Tx\| \| Sx\|  \leq \| T\| \| S\|.
\end{align*}
Thus $\| Tx\| = \| T\|$, $\| Sx\| = \| S\|$ and  $|[\langle Tx ,  Sx \rangle \xi, \xi]| = \| Tx\| \| Sx\|$. It follows from Theorem \ref{T1} that $T x\| Sx$.

$(iii) \Rightarrow(i)$: It is obvious.
\end{proof}
\end{proposition}

Popovici \cite[Theorem 2]{Dan} presented the following characterization.
\begin{lemma}\label{L4}
Let $\mathcal{F}$ be a Hilbert $\mathcal{A}$-module and $x, y \in \mathcal{F}$. Then the following statements are equivalent:
\begin{enumerate}
\item[(i)] $\| \theta_{x,x} +\theta_{y,y} \| = \| \theta_{x,x}\| + \| \theta_{y,y}\|$,
\item[(ii)] $\| \langle x, y \rangle\| = \| x\|\| y\|$.
\end{enumerate}
\end{lemma}

The following proposition is derived from Corollary \ref{Tete} and the above lemma.
 \begin{proposition}\label{P2}
Let $x, y \in \mathcal{E}$. Then the following statements are equivalent:
\begin{enumerate}
\item[(i)] $\theta_{x,x} \| \theta_{y,y}$,
\item[(ii)] there exists a sequence of basic vectors  $\{x_n\}^{\xi_n}$ in $\mathcal{E}$
 such that
\begin{align*}
\lim_n [\langle x_n, x \rangle \langle x, y \rangle \langle y , x_n \rangle \xi_n , \xi_n ] = \lambda \| x\|^2 \| y\|^2
\end{align*}
for some $\lambda \in \mathbb{T}$,
\item[(iii)] $\| \langle x, y \rangle\| = \| x\|\| y\|$.
\end{enumerate}
\begin{proof}

$(i)\Leftrightarrow(ii)$: The conclusion follows immediately from Corollary \ref{Tete}.

$(ii) \Rightarrow(iii)$: Suppose that $(ii)$ holds. we obtain
\begin{align*}
\| x\|^2 \| y\|^2 = \lim_n |[\langle x_n, x \rangle \langle x, y \rangle \langle y , x_n \rangle \xi_n , \xi_n ]| \leq \| x\| \|\langle x, y \rangle\| \| y\| \leq \| x\|^2 \| y\|^2,
\end{align*}
whence $\| \langle x, y \rangle\| = \| x\|\| y\|$.

$(iii) \Rightarrow(i)$. Let $\| \langle x, y \rangle\| = \| x\|\| y\|$. By $(ii) \Rightarrow(i)$ of Lemma \ref{L4}, we observe that $\| \theta_{x,x} + \theta_{y,y}\| = \| \theta_{x,x}\| + \| \theta_{y,y}\|$, which implies that $\theta_{x,x} \| \theta_{y,y}$.
\end{proof}
\end{proposition}

 \begin{corollary}\label{E21}
Let $x, y \in \mathcal{E}$. If $\langle x, y \rangle$ is a normal element, then the following statements are equivalent:
\begin{enumerate}
\item[(i)] $x\| y$,
\item[(ii)] $\| \langle x, y \rangle\| = \| x\|\| y\|$.
\item[(iii)] $\theta_{x,x} \| \theta_{y,y}$,
\item[(iv)] there exists a sequence of basic vectors  $\{x_n\}^{\xi_n},$ in $\mathcal{E}$
  such that
\begin{align*}
\lim_n [\langle x_n, x \rangle \langle x, y \rangle \langle y , x_n \rangle \xi_n , \xi_n ] = \lambda\| x\|^2 \| y\|^2
\end{align*}
for some $\lambda \in \mathbb{T}$,
\end{enumerate}
\begin{proof}
$(i)\Leftrightarrow(ii)$: It follows from the equivalence $(i)\Leftrightarrow(ii)$ of \cite[Theorem 2.6]{Z2}. Implications $(ii) \Leftrightarrow(iii)$ and $(iii)\Leftrightarrow(iv)$ follows from Proposition \ref{P2}.
\end{proof}
\end{corollary}

\begin{Example}\label{E2}
Let $\mathcal{H}$ be Hilbert space. Consider Hilbert $C^*$ -module $K(\mathcal{H})$ over itself via the inner product $\langle t, s \rangle = ts^*$. Fix $a \in K(\mathcal{H})$ and consider
the map $R_a: K(\mathcal{H}) \to K(\mathcal{H})$, defined by $R_a(t) = ta$, for every $ t \in K(\mathcal{H})$. Then $R_a$ is adjointable with the adjoint operator $R^*_{a}(s) = \langle s, a \rangle$. It is easy to see that $\theta_{t,s}= R_{t^*s}$. If $\langle t, s \rangle = ts^*$ is normal, then by Corollary \ref{E21},
we observe that $t\| s$ if and only if $R_{t^*t} \| R_{s^*s}$, and this holds if and only if there exists a sequence of basic vectors  $\{x_n\}^{\xi_n}$ in $K(\mathcal{H})$ such that $\displaystyle{\lim_n} [\langle x_n, t \rangle \langle t, s \rangle \langle s , x_n \rangle \xi_n , \xi_n ] = \lambda \| t\|^2 \| s\|^2$
and this occurs if and only if $\displaystyle{\lim_n }[|s|^2 v^*_{n} \xi_n , |t|^2 v^*_{n}\xi_n ] = \lambda\| t\|^2 \| s\|^2$ for some $\lambda \in \mathbb{T}$, or equvalently, $\|ts^* \| = \| t\|\| s\|$.
\end{Example}

As an immediate consequence of Proposition \ref{P2} and  Theorem 2.3 of \cite{AR3}, we present another characterization for pairs of compact operators in $\mathcal{L}(\mathcal{E}).$

 \begin{corollary}
 Let  $K(\mathcal{H})$ be a finite dimensional $C^*$-algebra and $\mathcal{E}$ be a full Hilbert $K(\mathcal{H})$-module
 and $x, y \in \mathcal{E}$ . Then the following statements are equivalent:
\begin{enumerate}
\item[(i)] $\theta_{x,x} \| \theta_{y,y}$,
\item[(ii)] there exists a basic vector  $\{v\}^{\xi}$ in $\mathcal{E}$
 such that
\begin{align*}
 [\langle v, x \rangle \langle x, y \rangle \langle y , v\rangle \xi, \xi ] =\lambda \| x\|^2 \| y\|^2
\end{align*}
for some $\lambda \in \mathbb{T}$,
\item[(iii)] $\| \langle x, y \rangle\| = \| x\|\| y\|$.
\end{enumerate}
\end{corollary}

\end{document}